\documentclass[preprintnumbers,amsmath,amssymb,twocolumn]{revtex4}
\usepackage[colorlinks=false,citecolor=red,filecolor=green,linkcolor=blue,pdfnewwindow=true]{hyperref}
\usepackage{hyperref}
\usepackage{titlesec}
\titleformat{\section}
	{\center\bfseries}{\makebox[30pt][l]{\thesection}}{0pt}{}
\titleformat{\subsection}
	{\flushleft\bfseries}{\makebox[30pt][l]{\thesubsection}}{0pt}{}
\titleformat{\subsubsection}
	{\flushleft\itshape}{\makebox[30pt][l]{\thesubsubsection}}{0pt}{}
\titlespacing\subsection{0pt}{12pt plus 4pt minus 2pt}{0pt plus 2pt minus 2pt}
\titlespacing\subsubsection{0pt}{12pt plus 4pt minus 2pt}{0pt plus 2pt minus 2pt}
\usepackage[font=small]{caption}
\usepackage{subcaption}
\usepackage{color}
\usepackage{graphicx}
\usepackage{makeidx}
\usepackage{bm}
\usepackage[title]{appendix} 
\usepackage{tabu}
\usepackage{float}
\usepackage{chemarr}
\usepackage{url}

\begin{document}

\title{Data driven approach to study the transition from dispersive to dissipative systems through dimensionality reduction techniques}
\author{Mairembam Kelvin Singh$^{1}$}
\author{A. Surjalal Sharma$^2$}
\author{N. Nimai Singh$^1$}
\author{Moirangthem Shubhakanta Singh$^1$}
\email{mshubhakanta@yahoo.com}
\affiliation{$^1$Department of Physics, Manipur University, Canchipur-795003, Imphal, Manipur, India.\\ $^2$Department of Astronomy, University of Maryland, College Park, MD 20742-2421, USA.}

\begin{abstract}
{\noindent}Complexity is often exhibited in dynamical systems, where certain parameters evolve with time in a strange and chaotic nature. These systems lack predictability and are common in the physical world. Dissipative systems are one of such systems where the volume of the phase space contracts with time. On the other hand, we employ dimensionality reduction techniques to study complicated and complex data, which are tough to analyse. The Principal Component Analysis (PCA) is a dimensionality reduction technique used as a means to study complex data. Through PCA, we studied the reduced dimensional features of the numerical data generated by a nonlinear partial differential equation called the Korteweg de Vries (KdV) equation, which is a nonlinear dispersive system, where solitary waves travel along a specific direction with finite amplitude. Dissipative nature, specific to that of the Lorenz system, were observed in the dimensionally reduced data, which implies a transition from a dispersive system to a dissipative system.\\ 

\noindent\textbf{Keywords:} Dispersive system; Dissipative system; Dimensionality reduction; Principal Component Analysis.
\end{abstract}
\maketitle

\section{Introduction}
{\noindent}Dynamical systems \cite{Birkhoff,Brin} that evolve with space and time show contrasting and varying features that uniquely define the system. These features are often found to be complex and cannot be studied easily. They are mathematically represented by linear or nonlinear differential equations and there are various analytical techniques to solve them. However, in certain cases analytical techniques are not feasible to produce a realistic result, wherein we use numerical techniques that provide us numerical data from which we can study the qualitative features of the system.\\

{\noindent}Dissipative \cite{Nicolis,Brogliato,France} and dispersive systems \cite{Whitham,Benjamin} are types of dynamical systems that exhibit complex behavior and rely on numerical techniques to produce realistic and reliable numerical solutions for better analysis. These two systems actually show contrasting features. The most common example of a dissipative system is the Lorenz strange attractor, which is based on a convection system in fluids with uniform depth and a linearly approximated temperature gradient between the layers of the fluid \cite{Saltzman,Lorenz}. The solutions of the Lorenz equations have irregular oscillations but always remain confined to a particular bounded region, giving a unique trajectory, called as ``strange attractor". For the case of dispersive system, we can cite the Korteweg de Vries (KdV) equation \cite{KdV} which describes the time evolution of dispersive waves. Its solutions known as ``solitons" belong to the class of solitary waves, that travel along a direction with a finite amplitude.\\

{\noindent}The numerical solutions of the above mentioned dynamical systems that show complex behavior are interpretable and can be derived using conventional numerical techniques \cite{Stuart} like the Runge-Kutta \cite{Bala,Jain,Shubha,Butcher} and multi-step methods \cite{Stuart}. However, analysis of large scale and complicated numerical data can be cumbersome. So, dimensionality reduction \cite{Reddy} and phase space reconstruction techniques \cite{King,Broomhead,Sharma} are also available for easier qualitative interpretation.\\

{\noindent}In \cite{Broomhead}, Takens' theorem \cite{Takens} and Singular Spectrum Analysis (SSA) \cite{Bertero,Pike} were used for phase space reconstruction of the Lorenz system. This technique was also used in \cite{Sharma} for studying the low dimensional behavior of the global magnetospheric dynamics. In this paper, we try to replicate the motive used in these past studies with a different technique called the Principal Component Analysis (PCA) \cite{Jolliffe,Bradde}. We also applied this technique to the dispersive system of the KdV equation and found that there is a transition to dissipative nature at reduced dimensions of the system. We solely focused on a data-driven approach by using the numerical solutions generated through conventional numerical techniques.\\

{\noindent}The paper is organised as follows. In section II, the general steps and procedure for implementing the technique that we used are described. Section III describes the various results that we obtained in two stages; the first part is to validate the technique by applying to the Lorenz system and obtaining similar results as done in previous work \cite{Broomhead}, the second part is to apply the technique to the dispersive KdV system. A brief discussion of the results that we obtained about the transition from the dispersive to dissipative system is also given.\\

\section{Methodology}
{\noindent}The first preliminary step is to produce the numerical solution of the system using a conventional numerical technique. While simple analytical technique like the separation of variables are more suited to linear PDEs, numerical techniques like the Runge-Kutta 4$^{\text{th}}$ order (RK4) method are more commonly used for nonlinear PDEs. Most complex systems have a nonlinear PDE or a coupled system of nonlinear PDEs associated with them. These PDE(s) can describe the time evolution of the dynamical variables involved. In some cases, the solution(s) of the PDE(s) describe a spatio-temporal evolution of the dynamical variables along space and time. In the two systems that we consider here, we applied the RK4 method \cite{Bala,Jain,Butcher} to generate the numerical solutions. This method is based on the general extrapolation equation, $Y_{i+1}=Y_i+mh$, where $m$ represents the slope and $h$ represents the step size. For a first order differential equation, $\displaystyle \frac{dY}{dX}=f(X,Y)$, with an initial point $(X_i, Y_i)$, the RK4 method uses the following formula for determining the numerical solution at a particular point:
\begin{eqnarray}
\label{RK4}
Y_{i+1}=Y_i+\left(\frac{m_1+2m_2+2m_3+m_4}{6}\right)h
\end{eqnarray}
where,
\begin{eqnarray}
\label{m}
\begin{split}
m_1&=f(X_i, Y_i)\\
m_2&=f\left(X_i+\frac{h}{2},Y_i+\frac{m_1h}{2}\right)\\
m_3&=f\left(X_i+\frac{h}{2},Y_i+\frac{m_2h}{2}\right)\\
m_4&=f\left(X_i+h,Y_i+m_3h\right)
\end{split}
\end{eqnarray}

{\noindent}However, for cases where there are higher order equations involved, we generally reduce them to a first order equation. After generating the numerical solution, we use it to construct the trajectory matrix \cite{Hassani,Darcy,Golyan} of the system. This step will depend on the system under consideration. For the KdV system, we treated the numerical solution directly as its trajectory matrix as the numerical solution is already a two dimensional data. Moreover, construction of trajectory matrix is commonly employed for one dimensional linear data.\\

{\noindent}After construction of the trajectory matrix, we then obtain the corresponding covariance matrix, which contains the covariances of all elements in the data \cite{Hui}. This is the first step in PCA and the step that differentiates the technique from other techniques like the SSA. We also use the Singular Value Decomposition (SVD) to enable us to classify the data according to the singular values of the covariance matrix, arranged in a decreasing order of their magnitude. This is followed by the extraction of the principal components, which we use to project the original data along the eigenvectors of the principal components. With this, we are able to obtain a dimensionally reduced data \cite{Li} and achieve to have a reconstructed phase space of the system.\\

{\noindent}The various steps involved in the technique we used are listed as follows:
\begin{enumerate}
\item Compute the numerical solution of the system using conventional numerical techniques. From the numerical data, obtain the corresponding trajectory matrix of the system, say $\bar{X}$.
\item Construct the covariance matrix of the numerical data. For any two variables $x$ and $y$, their covariance is given by $\sigma_{xy}^2=E[(x-\bar{x})(y-\bar{y})]$. For a single variable, its covariance is simply called variance, and is given by $\sigma_{x}^2=E[(x-\bar{x})^2]$ \cite{Hui}.
\item Apply Singular Value Decomposition (SVD) to the covariance matrix and the singular matrix so obtained will determine the covariances of the data arranged in order of their decreasing magnitude. This will enable us to extract the principal components with larger covariances.
\item Extract the principal components from the observations made using SVD. Then, construct the projection matrix, say $P$, using the eigenvectors corresponding to the singular values of the principal components \cite{Li}.
\item Use the projection matrix, $P$ to perform matrix dot product to the original matrix, $\bar{X}$ and obtain the dimensionally reduced matrix, say $X^\prime=\bar{X}P$, whose columns are the eigenvectors that contain the extracted features \cite{Li}.\\
\end{enumerate}

\subsection{Trajectory Matrix}
{\noindent}For the construction of the trajectory matrix, we take a linear time series $x(t)$ and transform it into a multi-dimensional series $X_0,...,X_K$, where each vector $X_i$ is called the $L$-lagged vector, the $L$ denoting the window length. $K$ will represent the number of columns of the trajectory matrix and $L=N-K+1$, where $N$ is simply the total number of data points or we can call it the number of `moments' in the time series. The window length $L$ is also such that $\displaystyle 2\le L \le \frac{N}{2}$.\\

{\noindent}Each $L$-lagged vector $X_i$ can be represented as $X_i=\left[x^{\prime}_i,...,x^{\prime}_{N-L+i} \right]$ and so the trajectory matrix will be given by \cite{Hassani, Darcy, Golyan},
\begin{eqnarray}
\label{traj_mat}
\bar{X}=\begin{bmatrix}
x^{\prime}_0 & x^{\prime}_1 & ... & x^{\prime}_{N-L}\\
x^{\prime}_1 & x^{\prime}_2 & ... & x^{\prime}_{N-L+1}\\
. & . & ... & .\\
x^{\prime}_{L-1} & x^{\prime}_L & ... & x^{\prime}_{N-1}\\
\end{bmatrix}
\end{eqnarray}

{\noindent}All the elements along the anti-diagonal are equal and so the trajectory matrix is a Hankel matrix \cite{Widom}.\\

\subsection{Principal Component Analysis}
{\noindent}PCA reduces the dimension of a data by classifying the various data points using their correlation or covariance. After the classification or ordering of the data from the highest variance to the lowest variance, we extract a few of the highly uncorrelated data points and project the original data along the dimensionally reduced dataset. This projection will allow us to obtain a new dataset free from noise and easier to analyse because of the reduced dimensions \cite{Hui, Li, Jolliffe,Bradde}.\\

{\noindent}PCA starts with the construction of the covariance matrix, which contains the covariances of all elements of a data \cite{Hui}. Those components with the larger covariances, which we can choose to use it for dimensionality reduction and projecting the original data along it, are called the principal components. The number of components that we choose, depends on the convenience of the data. For example, if we want to choose those components that captures 99$\%$ of the variance, then, we can use the condition, $\displaystyle \frac{\sum_{i=1}^k \sigma_i}{\sum_{i=1}^r \sigma_i}\ge 0.99$, where $k$ denotes the number of components that we are choosing, $\sigma_i$ denotes the $i^{\text{th}}$ singular value and $r$ denotes the total number of components \cite{Jolliffe}. Applying SVD on the covariance matrix and plotting the singular values will give the idea of how many singular values do we have to choose as principal components.\\

{\noindent}After choosing the $k$ principal components, we then construct a projection matrix  $P$, simply by arranging the eigenvectors corresponding to the singular values of the principal components. Finally, we can use this projection matrix to perform matrix dot product to the trajectory matrix $\bar{X}$ in \eqref{traj_mat} to obtain a dimensionally reduced matrix $X'=\bar{X}P$, whose columns are the eigenvectors of the reduced system that contains features of the entire system. Upon plotting the corresponding columns of the projected matrix pairwise, we can get the reconstructed phase space of the system.\\

\subsection{Singular Value Decomposition}
{\noindent}Singular Value Decomposition (SVD) is a technique to detect the part of a system with the most variation and use the subsequent part of the system to approximate the entire system with reduced dimensions by ordering the data points from the most variation to the least \cite{svd,Stewart}.\\

{\noindent}SVD begins with the fact that an $m\times n$ matrix $\psi$ can be decomposed into three matrices \cite{svd,Stewart}, two orthogonal matrices $L$, $R$ and a diagonal matrix $S$ such that,
\begin{eqnarray}
\label{SVD}
\psi_{mn} = L_{mm}S_{mn}R^T_{nn} 
\end{eqnarray}
where $L$ is the matrix containing the orthonormal eigenvectors of $\psi \psi^T$ and $R$ is the matrix containing the orthonormal eigenvectors of $\psi^T \psi$. The diagonal elements of the diagonal matrix $S$ are called the singular values, which are the square roots of the eigenvalues of $\psi^T \psi$ arranged in decreasing order. $S$ is called the singular matrix, $L$ and $R$ are respectively called the left singular and right singular matrices.\\

{\noindent}The idea behind SVD is that the singular values in $S$ and the orthonormal eigenvectors in $L$ and $R$ can be used to represent the original matrix $\psi$ as a linear combination of the three components. In this way, we can split a particular system and approximate it into lower dimensions from the analysis of the singular values so that we get a truncated version of the system.\\

\section{Results and Discussion}
{\noindent}Using the technique and steps described in the previous section, we implemented the same on two different systems. First, we implemented the technique on the dissipative Lorenz system. This part was intended to build the basis and define reference points for the other tasks. Second, we applied the technique to the dispersive KdV system. This was done to establish the transition study of dispersive to dissipative systems, purely based on a data driven approach. Both the systems considered are nonlinear and similar approach in applying the technique was used, with certain differences based on their numerical solutions.\\ 

\subsection{Lorenz System}

\begin{figure*}
\begin{center}
\label{Fig1}
\includegraphics[height=10.0cm,width=16.0cm]{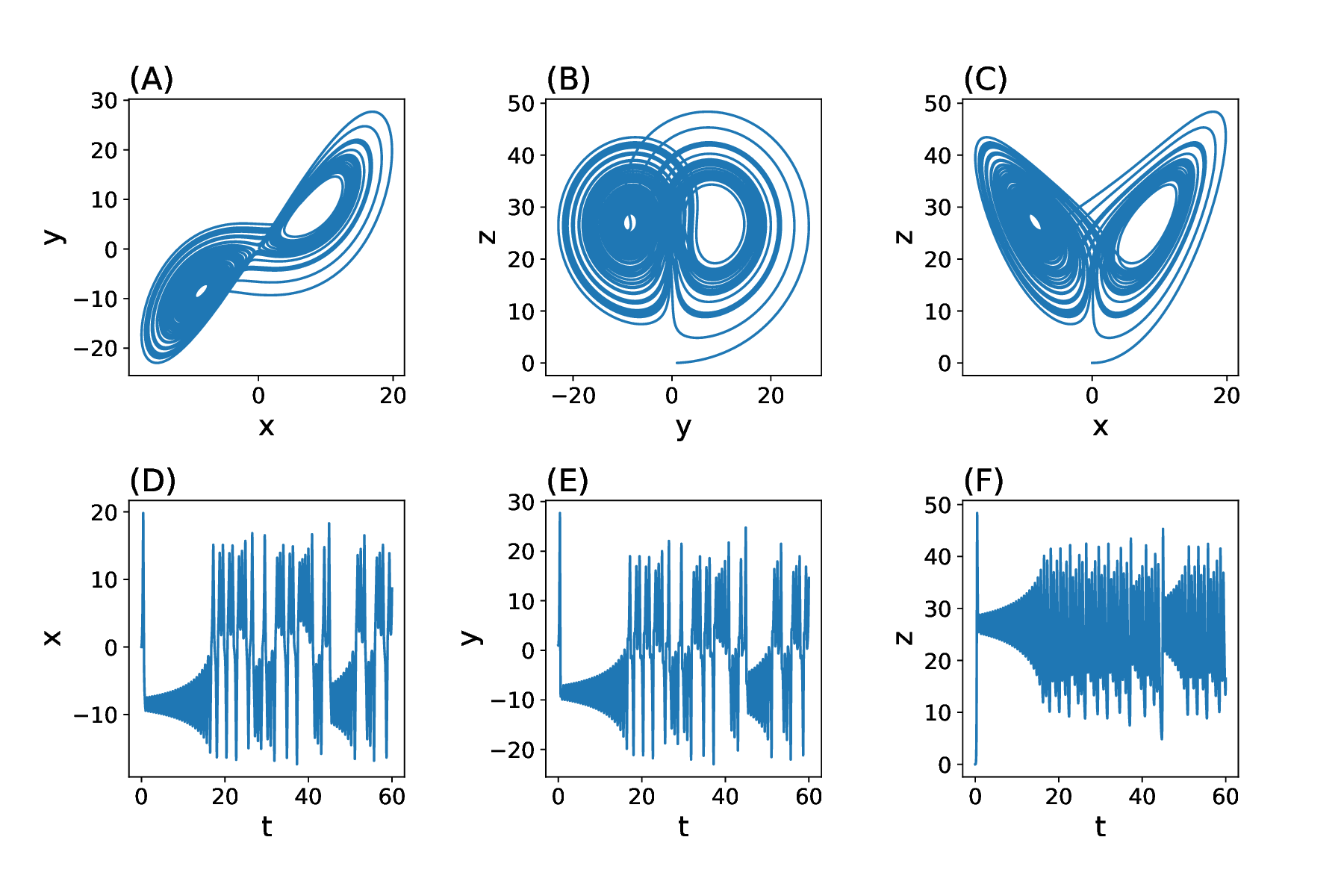}
\caption{(A)-(C): The phase diagrams of the three variables of Lorenz system, taken two at a time. We see that the trajectories confine to a bounded region. The points oscillate in an irregular and non-periodic manner. (D)-(F): These are the time series plots of each variable.}
\end{center}
\end{figure*}

{\noindent}Edward Lorenz, in 1963, simplified a model of convection system derived by Saltzman \cite{Saltzman}. This convection system was studied in a fluid system with uniform depth and a linearly approximated temperature gradient between the layers of the fluid. Lorenz derived the following set of equations \cite{Lorenz}:
\begin{eqnarray}
\label{Lorenz}
\begin{split}
\dot{x}&=\sigma(y-x)\\
\dot{y}&=rx-y-xz\\
\dot{z}&=xy-bz\\
\end{split}
\end{eqnarray}

{\noindent}Here, the variables $x$ is proportional to the intensity of the convective motion, $y$ is proportional to the temperature difference between the rising and declining fluid currents and $z$ is proportional to the deviation of the temperature gradient from linearity \cite{Lorenz}. $\sigma$ is called the Prandtl number, $r$ is called the Rayleigh number and $b$ gives the size of the region approximated by the system \cite{Record}. These three parameters are simply obtained as a result of dimensionless scaling to the convection equations derived by Saltzman \cite{Saltzman}. The dot simply denotes the first order derivative with respect to time.\\

{\noindent}Lorenz found out that the deterministic system in \eqref{Lorenz} had a strange dynamics. The solutions have irregular oscillations but always remain confined to a particular bounded region, giving a unique trajectory, called as ``strange attractor". This ``strange attractor" was later found out to be a fractal, having a fractional dimension between 2 and 3 \cite{Strogatz}.\\

{\noindent}Using the RK4 technique, the Lorenz equations \eqref{Lorenz} is solved with the initial values of the variables, $x_0=0$, $y_0=1$ and $z_0=0$; and the values of the parameters set at $\sigma=10$, $r=28$ and $\displaystyle b=\frac{8}{3}$. Consequently, we will get the time series data of the three variables, namely, $x(t)$, $y(t)$ and $z(t)$.\\

\begin{figure*}
\begin{center}
\label{Fig2}
\includegraphics[height=8.0cm,width=16.0cm]{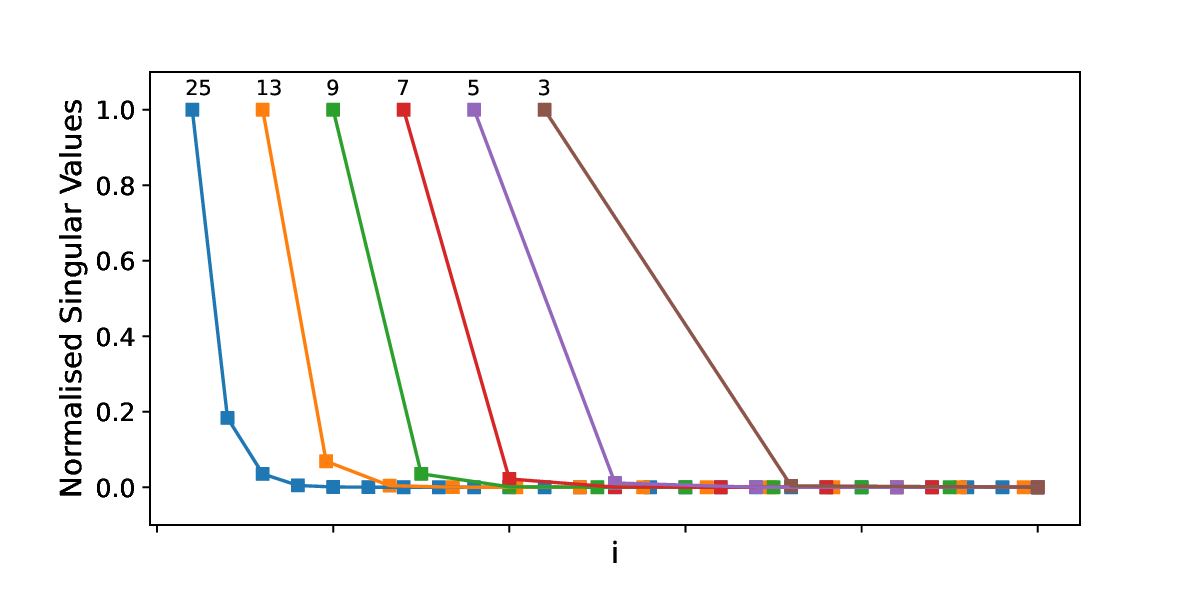}
\caption{Normalised singular values of the covariance matrix of the trajectory matrix of the time series data $x(t)$ of the Lorenz system. For each $\kappa$, there are two distinct parts, the deterministic part represented by the declining line and the noise floor represented by the horizontal and nearly flat part. As we decrease the value of $\kappa$, the number of singular values in the deterministic part decreases. The noise floor show the small differences between the smaller but finite singular values.}
\end{center}
\end{figure*}

{\noindent}The phase diagrams for the solutions are shown in Fig. 1(A)-(C). The time series plots of each variable are also shown in Fig. 1(D)-(F). Some distinctive properties derived from the figures are discussed as follows:
\begin{enumerate}
\item \textit{Non-periodic -} The trajectory of each variable does not approach a periodic or an equilibrium limit. There is an unstable equilibrium point at the centre of each loop around which the points rotates and then switches to the other unstable equilibrium point \cite{Record}.
\item \textit{Nonlinearity -} The Lorenz system \eqref{Lorenz} has two quadratic terms $xy$ and $xz$, which makes the system nonlinear \cite{Strogatz}.
\item \textit{Sensitive to Initial Conditions -} A characteristic feature of the Lorenz system is that a slight perturbation in the initial values of the variables can bring a drastic change to the entire outcome of the system's dynamics. As shown in Fig. 1(D)-(F), the trajectories which began with similar conditions, eventually depicts quite an unrelated behavior after some time. This property also makes long-term prediction impossible, unlike what classical physics does \cite{Record}.
\item \textit{Symmetry -} If we replace $(x, y)$ with $(-x, -y)$, there is no change in the equations. This shows that the Lorenz system is a symmetric system with symmetric solutions \cite{Strogatz}.
\item \textit{Dissipative -} Lorenz system is dissipative, which implies that volumes in the phase space contracts as time approaches infinity \cite{Strogatz}.
\item \textit{Fixed points -} For the Lorenz system \eqref{Lorenz}, $(0, 0, 0)$ is always a fixed point for all values of the different parameters. However, for $r>1$, there are two symmetric fixed points, $(\sqrt{b(r-1)},\sqrt{b(r-1)},r-1)$ and $(-\sqrt{b(r-1)},-\sqrt{b(r-1)},r-1)$. Also, for $r<1$, every trajectory approaches the origin $(0, 0, 0)$, as time goes to infinity. So, the origin is a globally stable fixed point for $r<1$ \cite{Strogatz}.\\
\end{enumerate}

{\noindent}The Lorenz system is referred to as ``chaotic" \cite{Chaos}, owing to the fact that for a certain range of parameters, there are no stable fixed points and no stable limit cycles, yet the trajectories confine to a bounded region \cite{Strogatz}. Among the values of the parameters that we used to solve the Lorenz system \eqref{Lorenz}, the value of $r$ is just past the Hopf bifurcation \cite{Hopf}, which implies a switch in the system's stability. Starting from the initial values, the solution sets into an irregular oscillation. The phase plot of $x$ vs $z$ looks like the wings of a butterfly \cite{Strogatz, Lorenz}. As we are viewing the trajectory from a two dimensional perspective, it seems like the trajectories cross each other. But, it won't be the case when viewed from three dimensions. The trajectory starts from the origin and spirals towards the central fixed point on the left. From there, it spirals outwards, and then jumps towards the right again after a certain period of time and then towards the left and so on indefinitely.\\

{\noindent}Using the time series data $x(t)$, we can construct the corresponding trajectory matrix $\bar{X}$. The number of samples considered in each window is varied as, $\kappa=(25,13,9,7,5,3)$ \cite{Broomhead}.\\

{\noindent}Next, we construct the covariance matrix of the trajectory matrix $\bar{X}$ and apply SVD to obtain the singular values, for each $\kappa$. The normalised singular values are plotted as shown in Figure 2. We see that for each $\kappa$, there are two distinct parts, the deterministic part represented by the declining line and the noise floor represented by the horizontal and nearly flat part. As we decrease the value of $\kappa$, the number of singular values in the deterministic part decreases. The noise floor show the small differences between the smaller but finite singular values. Practically, they should be zero but the presence of noise in the data makes them finite\cite{Sharma}.\\

\begin{figure*}
\begin{center}
\label{Fig3}
\includegraphics[height=8.0cm,width=18.0cm]{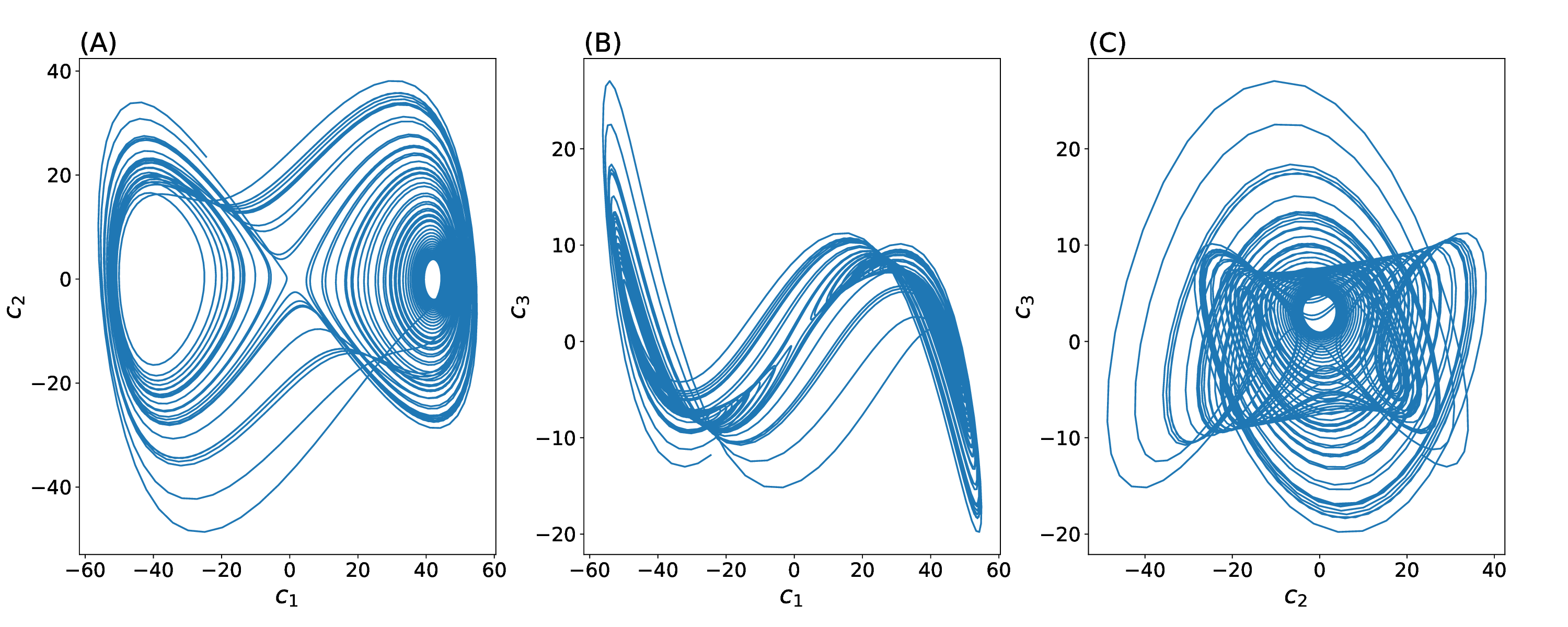}
\caption{Reconstructed phase space of the Lorenz system using three components of the reduced matrix $X^\prime$, namely, $c_1$, $c_2$ and $c_3$, plotted pairwise.}
\end{center}
\end{figure*}

{\noindent}To begin the reconstruction of the phase space and the subsequent dimensionality reduction, we shall use only the trajectory matrix corresponding to $K=25$. Then, we determine the singular values and eigenvectors of the covariance matrix and arrange the singular values according to their decreasing order. We select the top 3 singular values and their corresponding eigenvectors. These are called the principal components of the data and it captures most of the variance. Using these selected eigenvectors, we construct the projection matrix $P$. Finally, we can get the reduced matrix $X^\prime$ by performing dot product of the projection matrix $P$ and the trajectory matrix $\bar{X}$. To obtain the reconstructed phase space, we extract each column of the reduced matrix $X^\prime$ and plot them pairwise. The three columns or components of the reduced matrix are named as $c_1$, $c_2$ and $c_3$ \cite{Broomhead}. The pairwise projections of the reconstructed phase space are shown in Figure 3.\\

{\noindent}The reconstructed phase portraits obtained in Figure 3 are similar to that obtained in \cite{Broomhead}. This validates our technique and shows that we can conduct the process by using the covariance matrix, rather than by directly using the trajectory matrix. Now, after this validation we applied the technique to the nonlinear dispersive KdV system.\\

\subsection{Korteweg de Vries (KdV) System}
{\noindent}The Korteweg de Vries (KdV) equation \cite{KdV} is a non-linear 3$^{\text{rd}}$ order partial differential equation that describes the time evolution of dispersive waves. Its solutions known as ``solitons" belong to the class of solitary waves, that travel along a direction with a finite amplitude. These solitons can suffer collisions with one another and retain their shape and speed, showing their elastic nature \cite{Brauer, Zabusky}.\\

{\noindent}The KdV equation \cite{Drazin, Ablowitz} is given as,
\begin{eqnarray}
\label{KdV}
u_t + 6uu_x + u_{xxx} &=& 0 \nonumber\\
\Longrightarrow u_t + 3(u^2)_x + u_{xxx} &=& 0
\end{eqnarray}
Here,
\begin{eqnarray}
u &=& u(x,t);~~u_t = \frac{\partial u(x,t)}{\partial t}; \nonumber\\
u_x &=& \frac{\partial u(x,t)}{\partial x};~~u_{xxx}=\frac{\partial^3 u(x,t)}{\partial x^3}\nonumber\\
(u^2)_x &=& \frac{\partial}{\partial x}[u(x,t)]^2 \nonumber
\end{eqnarray}

{\noindent}The KdV equation can be numerically solved by imposing certain periodic boundary conditions. The Fourier pseudo-spectral method consisting of the continuous and discrete Fourier Transforms can be applied to such a periodic boundary problem \cite{Jie}. The numerical solution is based on the combination of the RK4 method \cite{Butcher,Bala,Jain} and Fast Fourier Transform (FFT) algorithm \cite{Jie, Trefethen, Dafermos}. To apply RK4, we need to convert the KdV equation to a first order ODE, which can be done after the discretization of the equation.\\

{\noindent}We set the initial condition $\displaystyle u(x, 0) = f(x) = \frac{v}{2}\text{sech}^2\left[\frac{\sqrt{v}}{2}x\right]$, which is based on the exact solution of the equation, $\displaystyle u(x, t) = \frac{v}{2}\text{sech}^2\left[\frac{\sqrt{v}}{2}(x-vt)\right]$, where $v$ is the velocity of the soliton and the periodic boundary condition $u(-l, t) = u(l, t)$. We, then, change the variable $x$ to $z = \frac{\pi x}{l}+\pi$ so that the domain is changed from $[-l, l]$ to $[0, 2\pi]$. This is needed to use the Discrete Fourier Transform (DFT) and apply the FFT algorithm. Making this change reduces the equation to,
\begin{eqnarray}
u_t + \frac{3\pi}{l}(u^2)_z + \frac{\pi^3}{l^3}u_{zzz} = 0 \nonumber
\end{eqnarray}

\begin{figure*}
\begin{center}
\label{Fig4}
     \centering
     \begin{subfigure}[b]{0.45\textwidth}
         \centering
         \includegraphics[width=\textwidth]{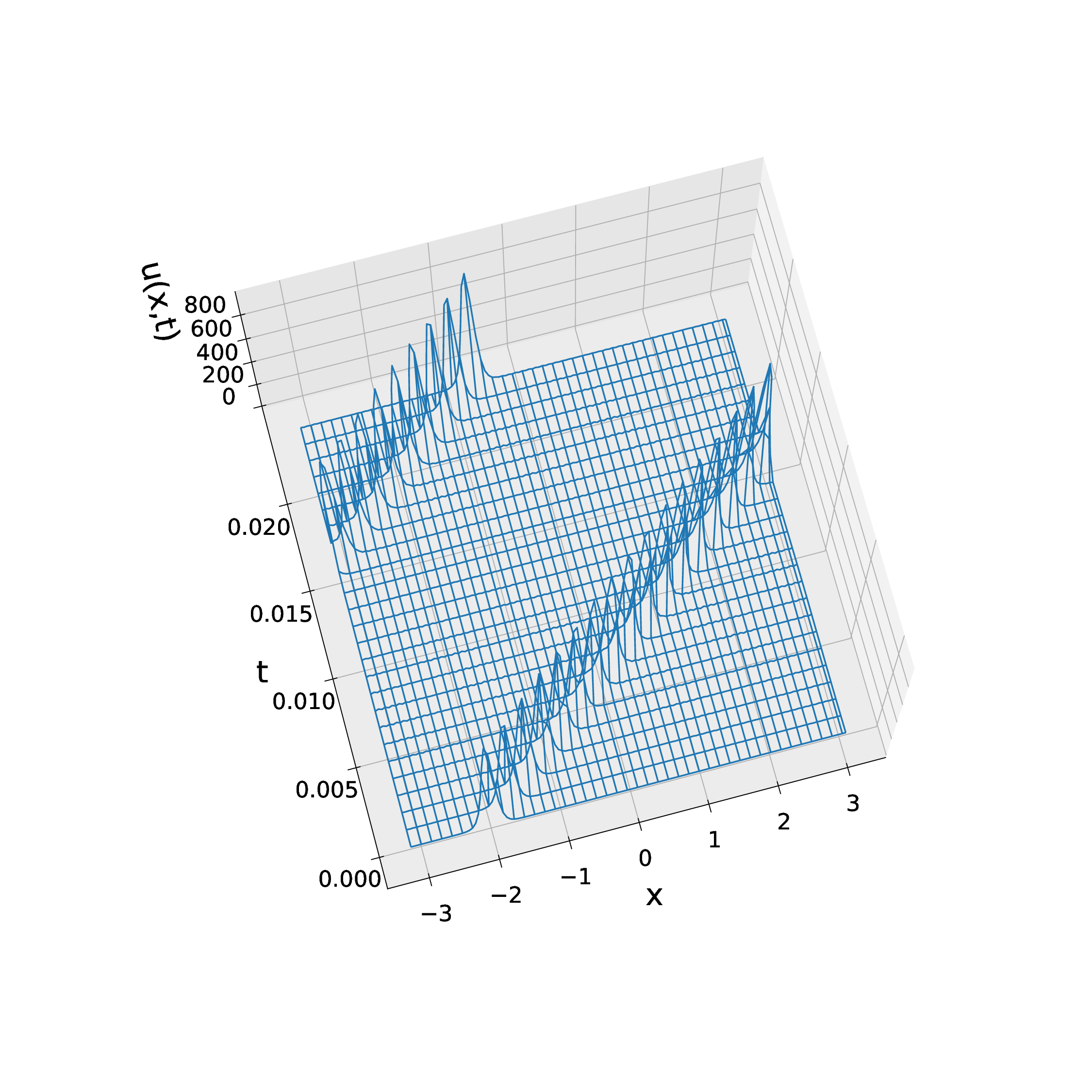}
         \caption{$N = 128 = 2^7$ and $v = 324$}
     \end{subfigure}
     \hfill
     \begin{subfigure}[b]{0.45\textwidth}
         \centering
         \includegraphics[width=\textwidth]{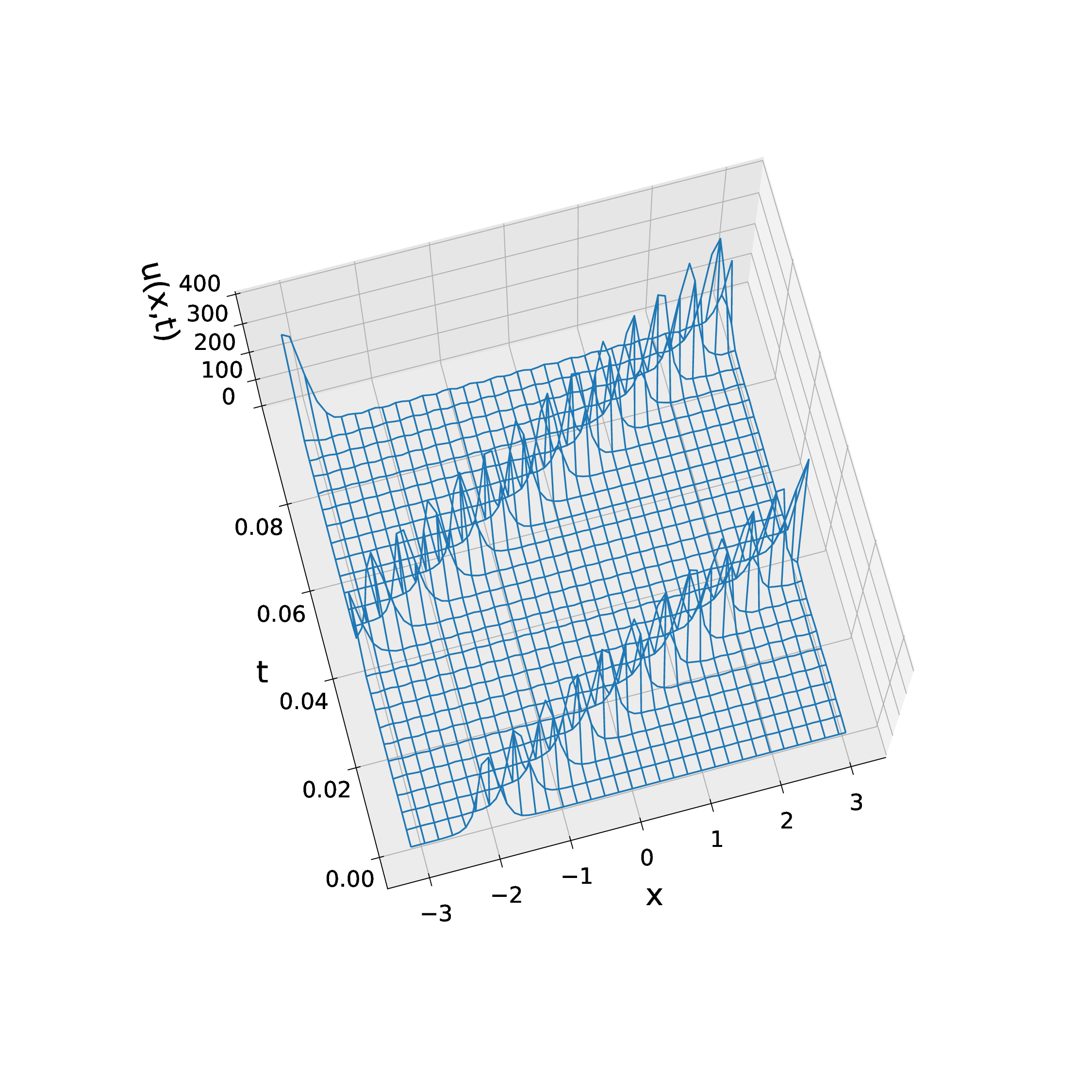}
         \caption{$N = 64 = 2^6$ and $v = 121$}
     \end{subfigure}
     \hfill
     \begin{subfigure}[b]{0.45\textwidth}
         \centering
         \includegraphics[width=\textwidth]{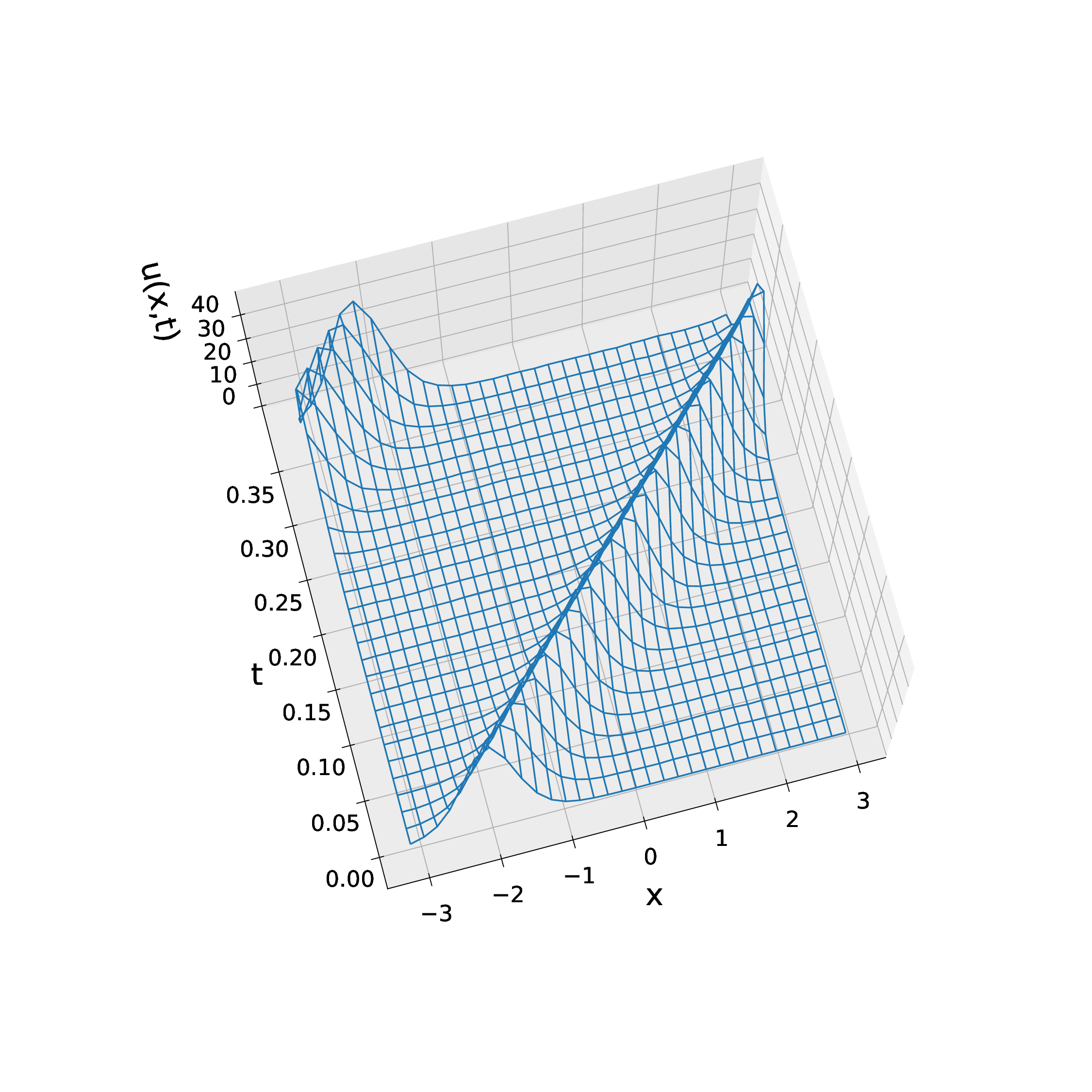}
         \caption{$N = 32 = 2^5$ and $v = 16$}
     \end{subfigure}
     \hfill
     \begin{subfigure}[b]{0.45\textwidth}
         \centering
         \includegraphics[width=\textwidth]{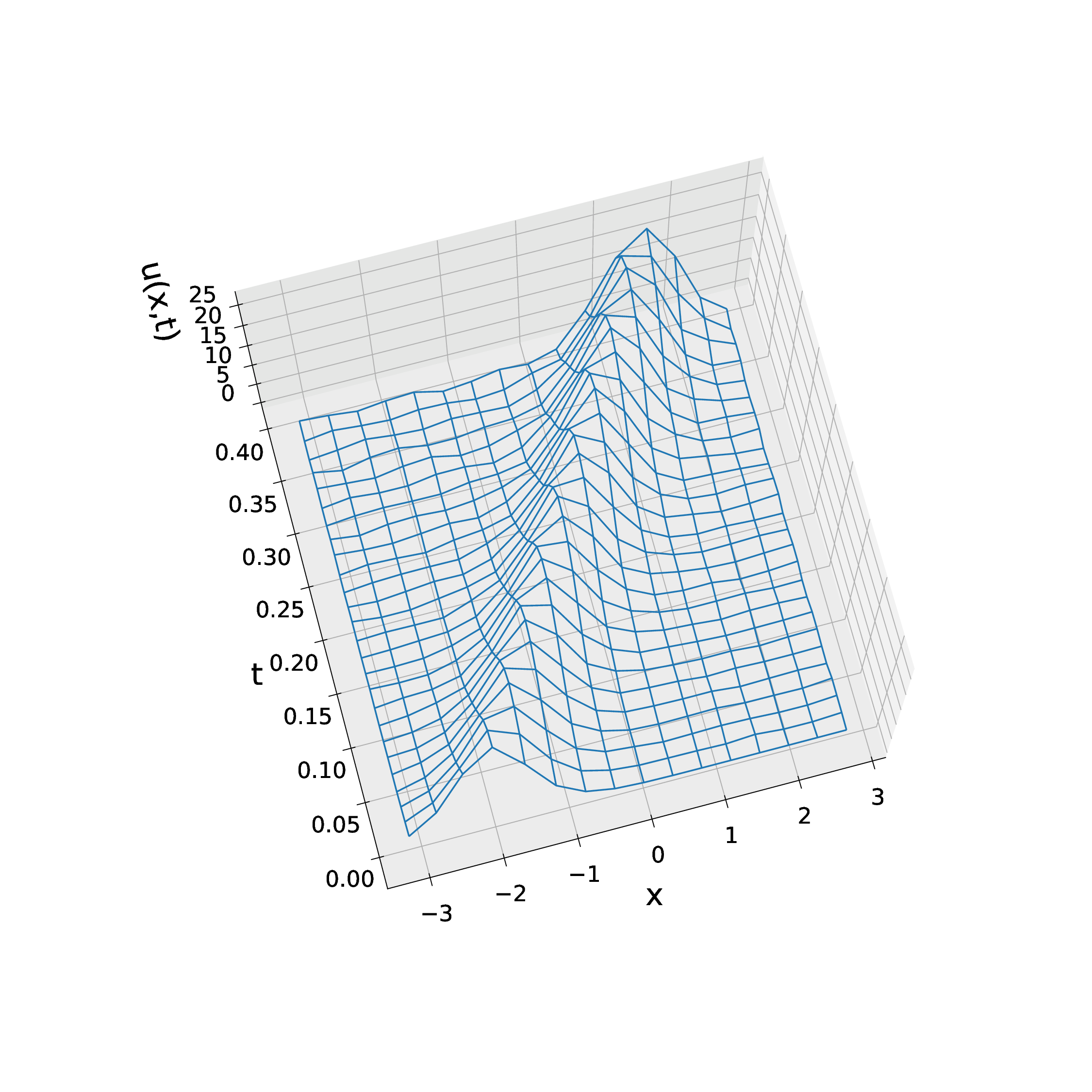}
         \caption{$N = 16 = 2^4$ and $v = 9$}
     \end{subfigure}
     \caption{Evolution of solitons with time and space for different values of $N$. We can observe the periodic behavior due to the periodic boundary conditions and the appearance of inherent noise indicated by the small spikes.}
\end{center}
\end{figure*}

{\noindent}Applying the DFT, we have, $\displaystyle u(x, t) \rightarrow u(z, t) = \sum_k \hat{w}_k (t) e^{ikz}$, where $\hat{w}_k$ is called the $k^{\text{th}}$ Fourier transform of $u(z, t)$. Then, using the relation, $\hat{w}_k^{(n)} = (ik)^n\hat{w}_k$, where the superscript ``$(n)$" on $\hat{w}_k$ denotes its $n^{\text{th}}$ order spatial derivative, we have,
\begin{eqnarray}
\frac{d\hat{w}_k(t)}{dt} + \frac{3i\pi k}{l}\widehat{(w^2)}_k - \frac{ik^3 \pi^3}{l^3}\hat{w}_k = 0 \nonumber
\end{eqnarray}

\begin{figure*}
\begin{center}
\label{Fig5}
\includegraphics[height=8.0cm,width=16.0cm]{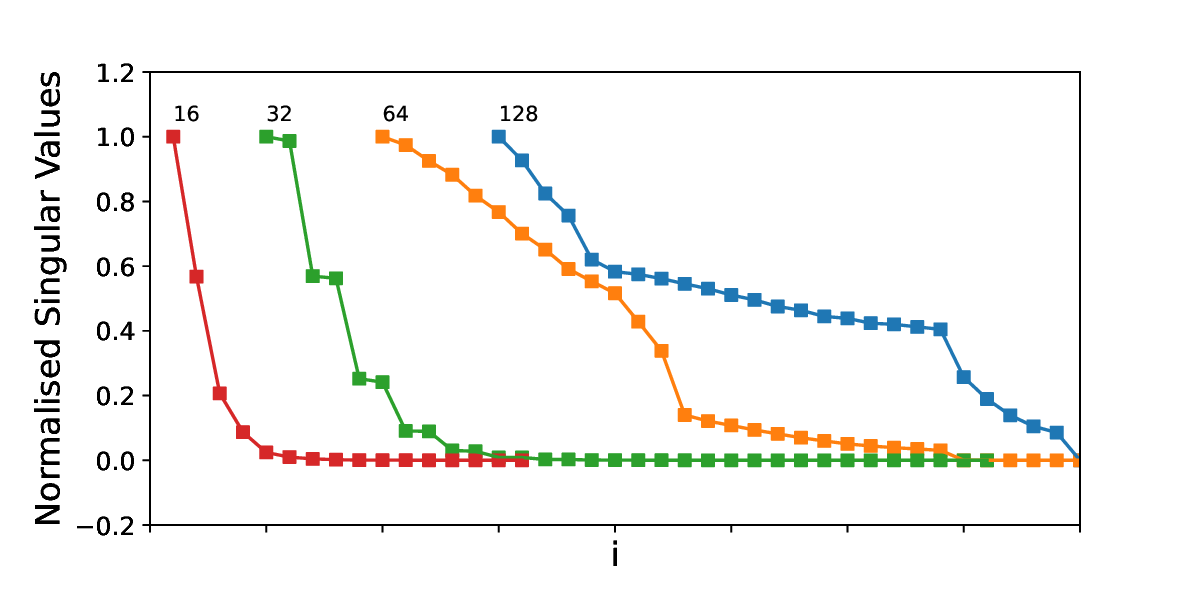}
\caption{Normalised singular values of the covariance matrix of the KdV trajectory matrix. For smaller values of $N$, we can see the larger covariance between the larger singular values. The noise floor is almost non-existent for larger values of $N$.}
\end{center}
\end{figure*}

{\noindent}Now, multiplying the above equation with the integrating factor $\displaystyle e^{-ik^3 \pi^3 t/l^3}$, we can reduce the equation as:
\begin{eqnarray}
\label{KdVRK4}
\frac{d}{dt}\left( e^{\frac{-ik^3 \pi^3 t}{l^3}} \hat{w}_k\right) = - \frac{3i\pi k}{l}e^{\frac{-ik^3 \pi^3 t}{l^3}}\widehat{(w^2)}_k \nonumber\\
\Longrightarrow \frac{d\hat{g}_k}{dt} = - \frac{3i\pi k}{l}\alpha F\left[\left\lbrace F^{-1}\left( \frac{\hat{g}_k}{\alpha}\right)\right\rbrace ^2\right]
\end{eqnarray}
where $\alpha = e^{\frac{-ik^3 \pi^3 t}{l^3}}$, $\hat{g}_k = \alpha \hat{w}_k$, $F[~~]$ denotes the Fourier transform and $F^{-1}(~~)$ denotes the inverse Fourier transform. Now, the RK4 method can be applied to \eqref{KdVRK4}. The corresponding computer code can be run using any programming language containing in-built modules to determine the FFT and inverse FFT (IFFT).\\

{\noindent}The particular FFT algorithm that we used for simulation is the ``decimation-in-time" algorithm for which the total length of the sequence must be a power of two \cite{Saidi}. So, the numerical simulation of the evolution of soliton(s) is made by initially setting up $N = 2^r$ grids for each particular value of time, where $r$ is an integer. These grids are known as Fourier collocation points \cite{Jie}. The step size for the values of time is set as $\displaystyle \Delta t = \frac{0.4}{N^2}$, so as to remain in the stability region \cite{Trefethen}.\\

{\noindent}For higher values of $N$, the velocity of the soliton is set at higher values because, if the velocity is small, we can almost never see the evolution with time and the soliton looks as if they are frozen. Also, the values are made a perfect square simply for smooth numerical calculation. The velocity of the soliton is reduced considerably for small number of grids because the results for FFT is too small for computation, at large velocities for small $N$. The step size and the upper limit of time is set accordingly so that each grid point is calculated using nearly 1000 iterations, for every value of $N$, except for the case of $N = 16$, where the number of iterations is reduced to 256 as we observe a very distorted numerical plot at larger time.\\

{\noindent}We begin with $N = 128 = 2^7$ grids, and then continue observing for smaller values, $N = 64 = 2^6$, $N = 32 = 2^5$ and $N = 16 = 2^4$. We observe the appearance of the periodic behavior of the wave due to the periodic boundary conditions and small spikes in the wave which indicates the presence of inherent noise in the system (Figure 4). From the results, we see a relationship between the appearance of periodicity and noise in the system for larger or smaller values of velocity of the soliton.\\

{\noindent}As the numerical data obtained from the numerical solutions are two dimensional, we treat it directly as the trajectory matrix of the data. For each plot corresponding to the values of $N$, we considered a numerical data consisting of only 26 time steps arranged along the row even though the whole simulation was run to approximately 1000 iterations. This is done mainly for two reasons, firstly, the time step was very small $(\approx 10^{-6}~\text{to}~10^{-3})$ and secondly to avoid complication. Thus the trajectory matrix will have 26 rows and $N$ number of columns.\\

\begin{figure*}
\label{Fig6}
     \centering
     \includegraphics[width=\textwidth]{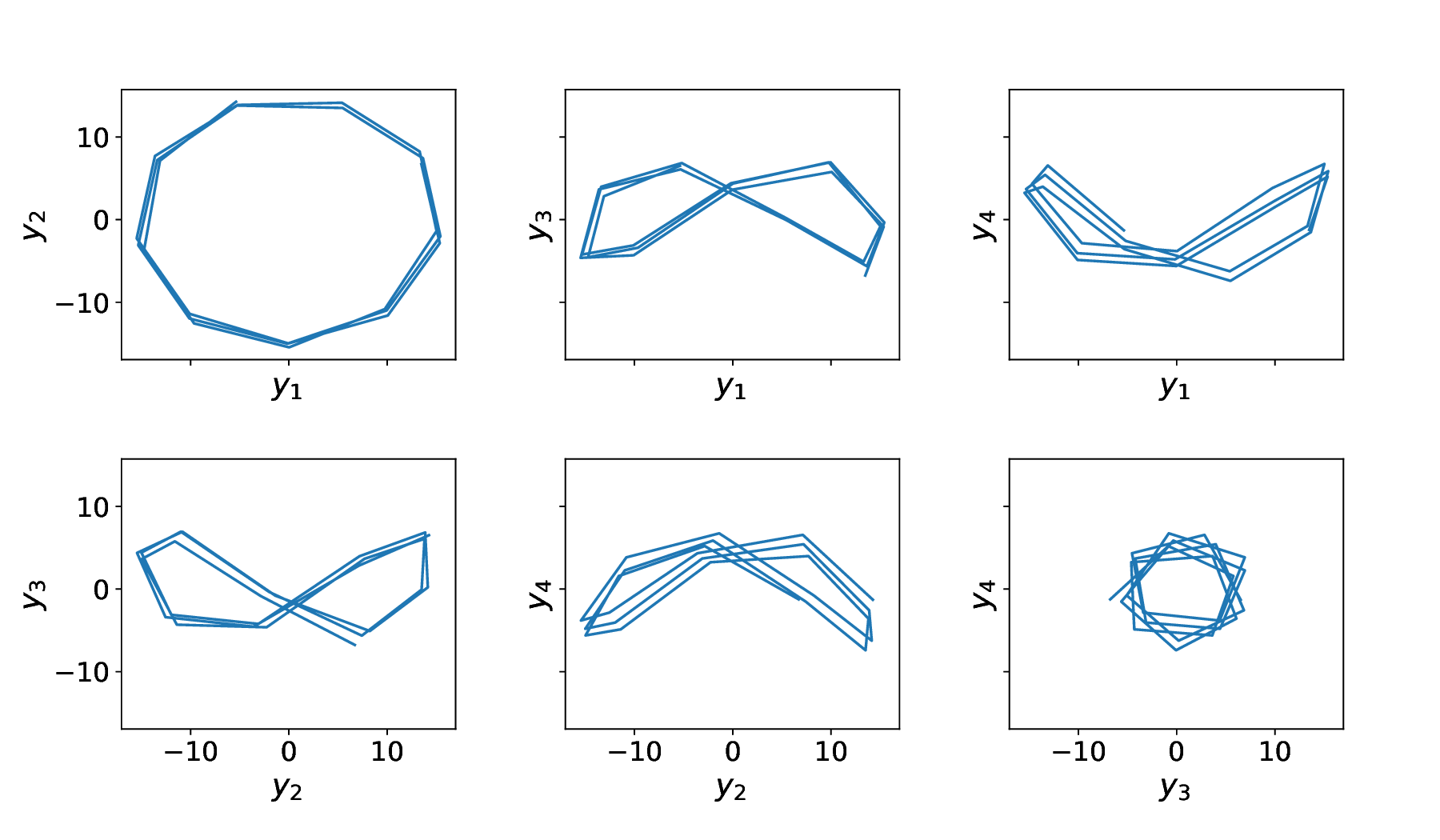}
     \caption{Reconstructed phase space of the reduced system of KdV Equations for $N = 16$. We can see the similarity of the phase space with that of a Lorenz system.}
\end{figure*}

\begin{figure*}
\label{Fig7}
     \centering
     \includegraphics[width=\textwidth]{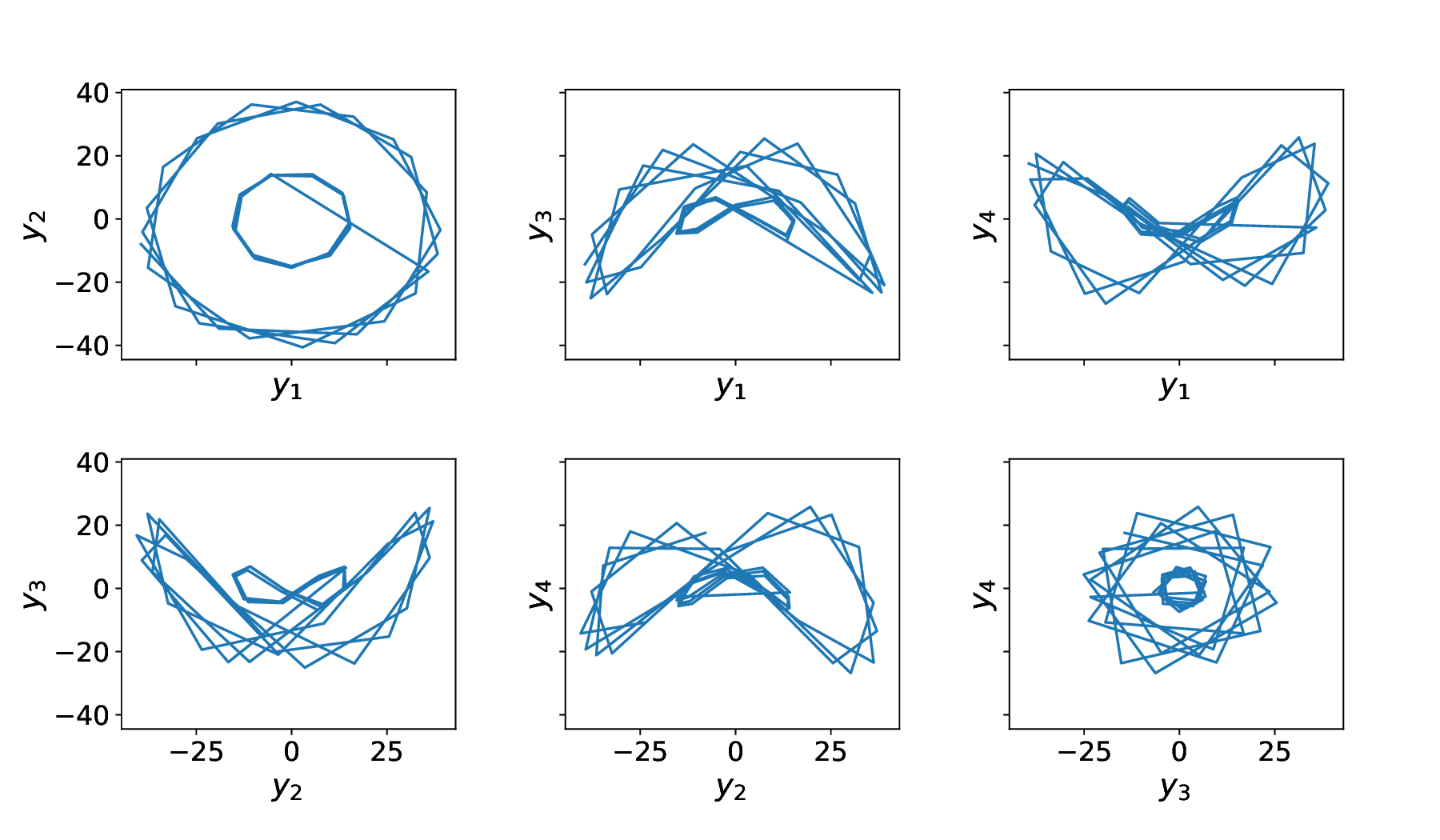}
     \caption{Reconstructed phase space of the reduced system of KdV Equations for $N = 32$. Here, a new feature in the form of a smaller closed loop is observed.}
\end{figure*}

\begin{figure*}
\label{Fig8}
     \centering
     \includegraphics[width=\textwidth]{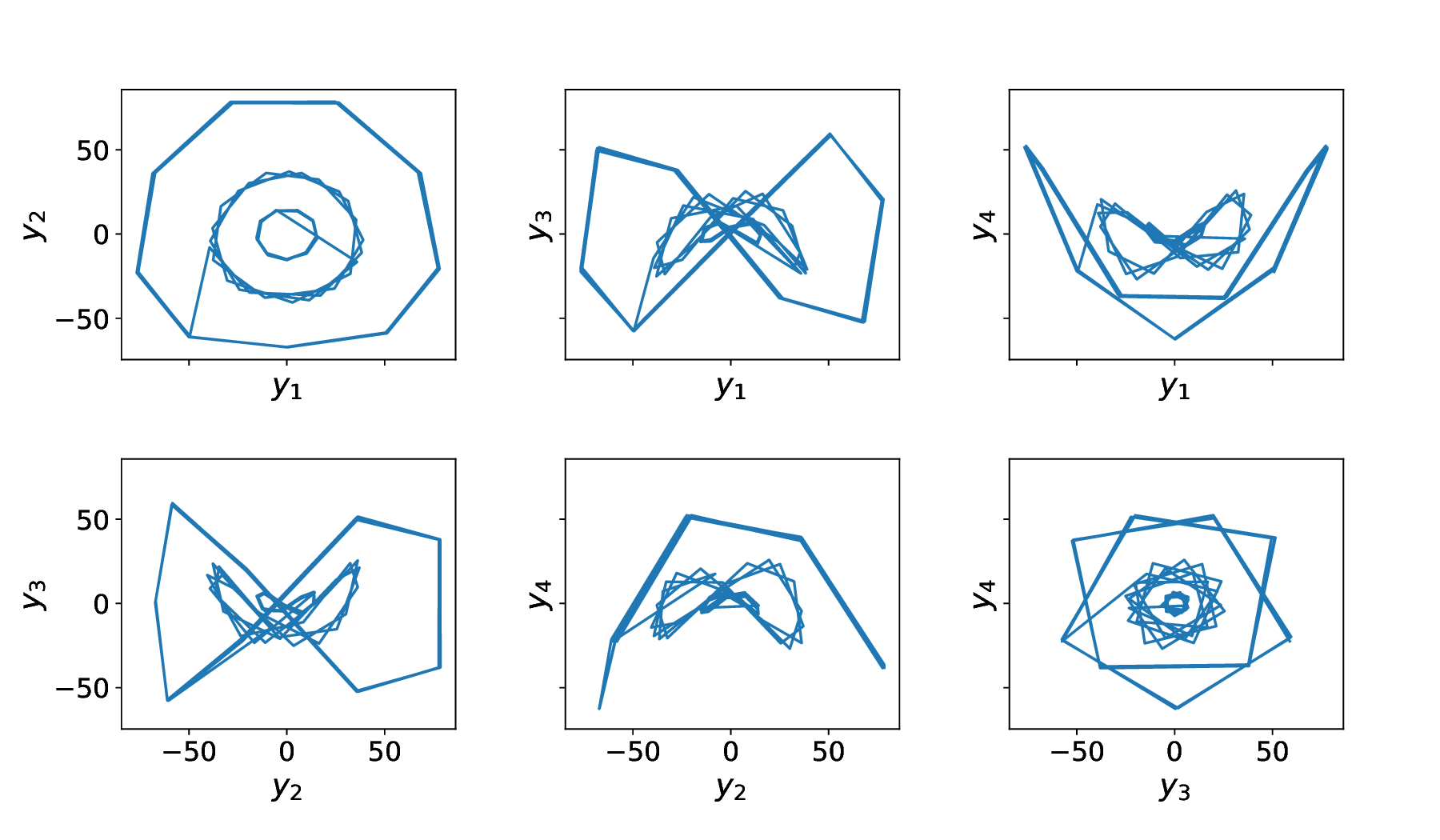}
     \caption{Reconstructed phase space of the reduced system of KdV Equations for $N = 64$. Here, we see the convergence of the same pattern.}
\end{figure*}

{\noindent}The covariance matrix is then constructed and SVD is applied. The plot of the normalised singular values is shown in Figure 5. We find a striking difference in the partition of the noise floor and the deterministic part for the Lorenz system as compared to the KdV system. In the Lorenz system, we see a steeper deterministic part and a flatter noise floor for larger samples while it is the opposite in the case of the KdV system. The noise floor is flatter and the deterministic part is steeper for smaller values of $N$ while this separation is quite non-existent for larger values of $N$.\\

{\noindent}For the KdV system, we choose four principal components and the original data can be projected along the eigenvectors of the projection matrix. The reduced matrix will have four columns or components and they are named as $(y_1, y_2, y_3, y_4)$. We stick to the cases of $N = 16, 32, 64$ and the pairwise plots of these components are made. These are shown in Figures 6, 7 and 8 respectively for each value of $N$ that we considered.\\

{\noindent}The reconstructed phase space of the KdV system resembles closely to that of a Lorenz system. The phase portraits also show closed concentric loops as we increase the values of $N$. Specifically, for $N=16$, there is only one loop, for $N=32$, there are two loops and for $N=64$, there are three loops. In making these phase space reconstruction, we decreased the velocities of the solitons from the previously set values. In doing so, we were able to observe a similarity in the pattern with that of a Lorenz system and the convergence of the similar pattern as $N$ is increased. The previous velocities of solitons used were 9, 16 and 121, respectively for $N = 16, 32, 64$ and they are now changed to 4, 9 and 16, respectively. The previous values were set so as to view a distinctive pattern of evolution of the solitons while working on the numerical solutions. As the amplitudes of the solitons depend directly on their velocities, decreasing the velocities will obviously decrease the amplitudes and we would have not seen the evolution that we saw in the plots of Figure 4.\\

{\noindent}The decrease in the velocities does not necessarily mean a change in the behavior of evolution of the solitons. It simply means that the amplitudes of the solitons are decreased considerably. There won't be any affect in the principle of the technique we were applying by the change of this parameter. This is in fact supported by the fact that the singular values of the covariance matrix show almost a similar pattern when the velocities are changed. The only distinctive change that can be recognised is the fact that the curves for $N = 64~\text{and}~128$ looks more definitive with respect to the separation of the deterministic part and the noise floor (Figure 9).\\

\begin{figure*}
\label{Fig9}
     \centering
     \includegraphics[width=\textwidth]{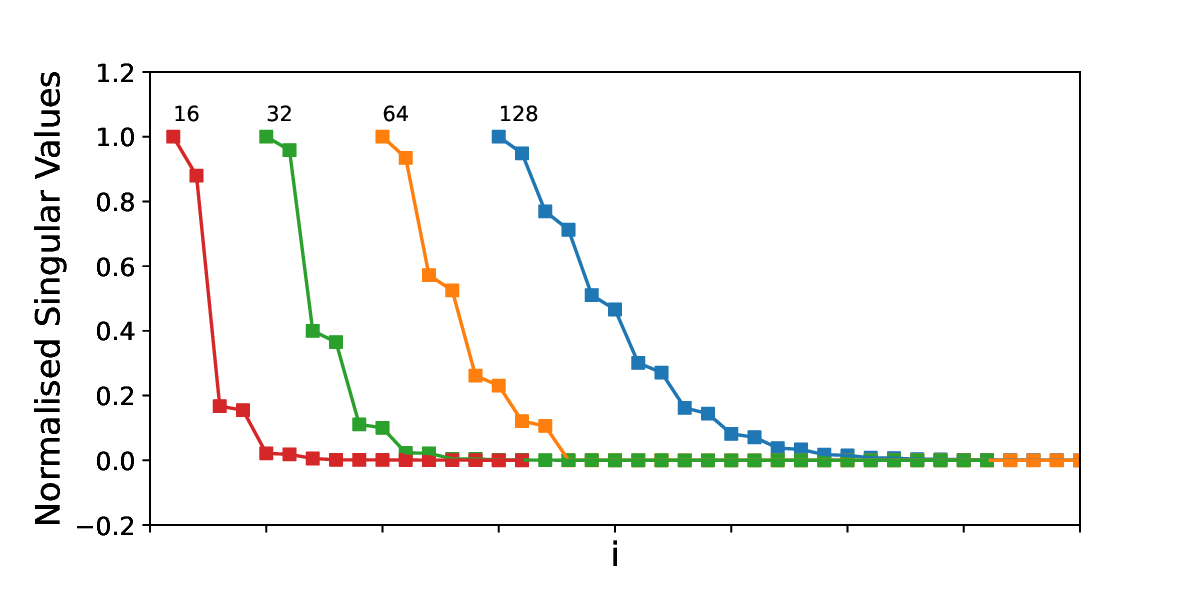}
     \caption{Normalised singular values of the covariance matrix of the KdV trajectory matrix for small velocities and hence small amplitudes of the solitons. The curves for $N = 64~\text{and}~128$ looks more definitive with respect to the separation of the deterministic part and the noise floor, as compared to that of Figure 5. }
\end{figure*}

{\noindent}So, we can conclude that at very small velocities and hence small amplitudes of the solitons, the dispersive KdV system transitions to a dissipative system. This is evident from the fact that the reconstructed phase portraits resembles that of the Lorenz system. Also, for higher values of $N$, we see concentric closed loops, which signifies the nature of volume contraction of a dissipative system.\\

\section{Conclusion}
{\noindent}The technique of Principal Component Analysis (PCA) is found to be effective in studying the reduced dimensional features of dynamical systems. The reconstruction of phase space of a dissipative system like the Lorenz system shows a similar result to what researchers have obtained in the past. However, when we applied the same technique to a dispersive system like the KdV system, we can observe a transition to a dissipative nature.

{\noindent}The reconstructed phase space of the KdV system using the principal components show a similar pattern to that of a Lorenz system. The trajectory of the phase space does not approach a stable equilibrium point, but remains confined to a specific boundary. As we increase the grid size $(N)$, we also observed that the pattern shows concentric closed loops, with the trajectory jumping from one loop to an inner loop, signifying a contraction in volume as the trajectory evolves. These are the characteristic features mentioned in the discussion of the Lorenz system.

{\noindent}However, the above observations are made by switching the values of velocities of the solitons. The values of velocities used for producing the numerical solutions were set to facilitate smooth numerical computations. These values could not produce the transitional characteristics of the reconstructed phase space. So, to display a good view of the concentric and convergent pattern, we decreased the values of the velocities. Quantitatively, this will only decrease the amplitudes of the solitons but will have no affect on the qualitative behavior of the spatio-temporal evolution of the solitons.\

{\noindent}Thus, we can conclude that the reduced dimensional features of a dispersive system shows dissipative behavior. In other words, our study reveals the dissipative nature of a dispersive system embedded in reduced or lower dimensions. This transition is also possible only at small amplitudes of the solitons. We completely emphasised on a data driven approach based on the numerical solutions provided by the conventional numerical techniques.\\

{\noindent}{\bf Acknowledgements} \\
{\noindent}M.K. Singh is a Junior Research Fellow (JRF) under the National Fellowship for Scheduled Castes Students (NFSC) scheme of University Grants Commission (UGC), Delhi, India. The simulations were done through Python 3 programs written using the Spyder environment of Anaconda distribution. Simulations that require large computational tasks were done in HPE ProLiant DL380 Gen10 Plus installed in the Department of Physics, Manipur University.\\
{\noindent}A.S. Sharma passed away before the completion of this manuscript. His guidance and motivation throughout the development of this paper is immensely acknowledged.\\

{\noindent}{\bf Author Contributions:}\\
{\noindent}A.S. Sharma conceptualized the work. M.K. Singh and M.S. Singh did the analytical, computational work and preparation of the figures. M.K. Singh, N.N. Singh and M.S. Singh wrote the manuscript, read, analysed the results, and approved the final manuscript.\\

{\noindent}{\bf Additional Information} \\
\textbf{Competing interests:} \\ The authors declare no competing interests.

\end{document}